\documentclass[12pt,twoside]{amsart}
 
\usepackage{amsmath} 
\usepackage{amssymb} 
\newtheorem{thm}{Theorem} 
 
\newtheorem{lemma}[thm]{Lemma}

\newtheorem{preremark}[thm]{Remark} 
 
\def\N{{\mathbb N}}

\def\R{{\mathbb R}} 
\def\a{{\alpha}} 
 
\newcommand{\CVD}{\hfill $\rule{2.5mm}{2.5mm}$} 
\newcommand{\PF}{\noindent{\bf Proof. }} 
 
\title[Fractional Laplacian phase transitions]{ 
Fractional Laplacian phase transitions\\
and boundary reactions:\\ 
a geometric inequality\\ 
and a symmetry result} 
 
\author[Y. Sire]{Yannick Sire} 
\author[E. Valdinoci]{Enrico Valdinoci}

\begin{document} 
\begin{abstract} 
We deal with symmetry
properties for solutions
of nonlocal equations of the type  
\begin{equation*} 
(-\Delta)^s v= f(v)\qquad 
{\mbox{ in $\R^n$,}} 
\end{equation*} 
where $s \in (0,1)$ and
the operator $(-\Delta)^s$ is the so-called fractional 
Laplacian. 
%% It is a pseudo-differential operator with  
%% symbol~$|\eta|^{2s}$. 
  
The study of this 
nonlocal equation is made via a careful 
analysis of the following degenerate elliptic equation   
$$
\left\{ 
\begin{matrix} 
-{\rm div}\, (x^\a \nabla u)=0 \qquad 
{\mbox{ on $\R^n\times(0,+\infty)$}} 
\\ 
-x^\a u_x = f(u) 
\qquad{\mbox{ on $\R^n\times\{0\}$}}\end{matrix} 
\right.$$
where $\a \in (-1,1)$.  
 
This equation is related to the fractional 
Laplacian since the  
Dirichlet-to-Neumann operator~$\Gamma_\a:  
u|_{\partial \R^{n+1}_+} \mapsto  
-x^\a u_x |_{\partial \R^{n+1}_+} $ 
is 
$(-\Delta)^{\frac{1-\a}{2}} $.

More generally, we study the so-called boundary reaction equations 
given by  
 \begin{equation*}\left\{ 
\begin{matrix} 
-{\rm div}\, (\mu(x) \nabla u)+g(x,u)=0 \qquad 
{\mbox{ on $\R^n\times(0,+\infty)$}} 
\\
-\mu(x) u_x = f(u) 
\qquad{\mbox{ on $\R^n\times\{0\}$}}\end{matrix} 
\right.\end{equation*}   
under some natural assumptions on the diffusion coefficient
$\mu$ and on
the nonlinearities $f$ and $g$. 

We prove a geometric formula of 
Poincar\'e-type for stable solutions, from which we 
derive a symmetry result
in the spirit of a conjecture of De Giorgi. 
\end{abstract} 
 
\maketitle 
\tableofcontents 
 
\bigskip\bigskip 
 
\bigskip 
 
\bigskip\bigskip 
 
\section*{Introduction} 
 
This paper
is devoted to some geometric results on the following equation  
\begin{equation}\label{lapFrac} (-\Delta)^s v= f(v)\qquad {\mbox{ in  
$\R^n$.}} \end{equation} 

The operator $(-\Delta)^s$ is the  
fractional  
Laplacian and it 
is a pseudo-differential operator with symbol~$|\eta|^{2s}$, with~$s\in 
(0,1)$ -- here, $\eta$ denotes the variable in the frequency 
space. This operator, which is a nonlocal operator, can also be defined, 
up to a multiplicative constant, by the  
formula 
\begin{equation}\label{kernel} 
(-\Delta)^s v (x)=P.V. \int_{\R^n}  
\frac{v(x)-v(y)}{|x-y|^{n+2s}}\,dy, 
\end{equation} 
where $P.V.$ stands  
for  
the Cauchy principal value (see \cite{landkof}
for further details).

Seen as an operator acting on distributional spaces, the
quantity $(-\Delta)^s v$ is well-defined as long as $v$ belongs to the
space 
\begin{equation*}
\mathcal{L}_s=\left \{v \in
  \mathcal{S}'(\R^n),\,\,\int_{\R^n}
  \frac{|v(x)|}{(1+|x|)^{n+2s}}\,dx<\infty\right \} \bigcap
C^2_{\rm{loc}}(\R^n).
\end{equation*}
Notice in particular that smooth
bounded functions are admissible for the
fractional Laplacian. The $L^1$ assumption allows to make the integral
in \eqref{kernel} convergent at infinity, whereas the additional
assumption of $C^2_{\rm loc}$-regularity is here to make sense to the
principal value\footnote{For $v\in C^2_{\rm loc}(\R^n)$,
the singular integral in \eqref{kernel} makes sense for any~$s
\in(0,1)$.
Of course, it is possible to weaken such assumption
depending on the values of $s$.} 
near the singularity.
 
{F}rom a probabilistic point of view, the  
fractional Laplacian is the infinitesimal 
generator of a Levy process (see, e.g.,~\cite{B}). 
 
This type of diffusion operators arise in several areas such as  
optimization~\cite{DL}, flame propagation~\cite{CRS} 
and finance~\cite{CT}. Phase transitions driven by
fractional Laplacian-type boundary effects have also been
considered in \cite{ABS, Mar} in the Gamma convergence
framework. Power-like nonlinearities for boundary reactions
have also been studied in \cite{POW}
 
In this paper, 
we focus on an analogue of the De Giorgi conjecture~\cite{DeG} 
for  
equations of the  
type \eqref{lapFrac}, 
namely, whether or not ``typical'' solutions 
possess one-dimensional symmetry. 
 
One of the main  
difficulty of the analysis of 
this operator is its nonlocal character.  
However, it is a well-known fact 
in harmonic analysis that the power $1/2$  
of the Laplacian is the boundary operator of 
harmonic functions in the half-space.  
 
In~\cite{cafS}, the equivalence
between \eqref{lapFrac} and the $\a$-harmonic extension in the half-space 
has recently been proved.
More  
precisely, if one considers the boundary reaction problem    
\begin{equation}\label{bdyFrac} 
\left \{
\begin{matrix} 
{\rm div}\, (x^\a \nabla u)=0 \qquad 
{\mbox{ on $\R^{n+1}_+ 
:=\R^n\times(0,+\infty)$}} 
\\
-x^\a u_x = f(u)  
\qquad{\mbox{ on $\R^n\times\{0\}$,}}\end{matrix}\right . \end{equation} 
it is proved in~\cite{cafS} 
that, up to a normalizing factor,
the Dirichlet-to-Neumann operator  
$\Gamma_\a:  
u|_{\partial \R^{n+1}_+} \mapsto  
-x^\a u_x|_{\partial \R^{n+1}_+} $ 
is precisely $(-\Delta)^{\frac{1-\a}{2}} $ and then that $u(0,y)$ is a
solution of 
$$(-\Delta)^\a u(0,y)=f(u(0,y)). $$ 
 
Note that 
the condition~$\frac{1-\a}{2}=s\in(0,1)$ in~\eqref{lapFrac} 
reduces to~$\alpha\in(-1,1)$. 
 
Qualitatively, the result of~\cite{cafS} 
states that one can localize the  
fractional  
Laplacian by adding an additional variable. 
This argument plays, for instance, a crucial  
role in the proof of full regularity of the 
solutions of the quasigeostrophic model as given by~\cite{CV}
and in the free boundary analysis
in~\cite{CSS}.
 
The  
operator ${\rm div}\,( x^\a \nabla )$ is elliptic degenerate. 
However, since $\a \in (-1,1)$, the weight $x^\a$ is integrable at $0$. 
This type of weights falls into the category of $A_2$-Muckenhoupt  
weights (see, for instance,~\cite{muck}), 
and an almost complete theory for these 
equations is available 
(see \cite{FKS, FJK}). In particular, 
one can obtain H\"older regularity,  
Poincar\'e-Sobolev-type estimates, Harnack and boundary Harnack  
principles.        

In the present paper, we 
want to give a geometric insight of the phase  
transitions for equation~\eqref{lapFrac}. Our goal 
is to give a geometric  
proof of the one-dimensional symmetry result for 
fractional boundary reactions in dimension~$n=2$, 
inspired by 
De Giorgi conjecture and 
in the spirit of the proof of 
Bernstein Theorem given in~\cite{giusti}. 
 
A similar De Giorgi-type result 
for boundary reaction 
in dimension~$n=2$  
has been proven in \cite{CSM}  
for~$\a=0$,
which corresponds 
to the square root of 
the Laplacian in~\eqref{lapFrac}.
The 
technique of \cite{CSM} will be adapted to 
the case~$\a\ne0$ in the forthcoming~\cite{CS}. 
 
However,  
the proofs in~\cite{CSM, CS} 
are based on different methods (namely, a Liouville-type 
result inspired 
by~\cite{BCN, AC, AAC} and a  
careful analysis of the  
linearized equation around a solution) 
and they are quite technical and long. 
Our techniques also gives some geometric
insight on more general types of boundary reactions
(see equation \eqref{eq1-provv} below).
\bigskip

There has been a large number of works   
devoted to the symmetry properties of 
semilinear equations for the standard Laplacian.  
In particular, De Giorgi conjecture on the flatness of level sets 
of standard phase transitions  
has been studied in low dimensions in \cite{AAC,AC,BCN,GG1,GG2}. 
The conjecture has also been settled in \cite{savin} 
up to dimension $8$ under an additional assumption 
on the profiles at infinity. 
 
Here we give a proof  
of analogous symmetry properties 
for  
phase transitions driven by fractional 
Laplacian as in \eqref{lapFrac}. 
Such proof 
will be rather simple and short, 
with minimal assumptions (even on the nonlinearity~$f$ 
which can be taken here to be just locally Lipschitz) 
and it 
reveals some geometric aspects of the equation. 
 
Indeed, our proof, which is based on the recent work~\cite{FSV},
relies heavily on a Poincar\'e-type inequality which involves the  
geometry  
of the level sets of $u$. 
\bigskip 

Most of our paper will focus
on the boundary reaction equation in \eqref{bdyFrac}
(and, in fact, on the more general framework of
\eqref{eq1-provv} below). We recall that \eqref{bdyFrac}
still exhibits
nonlocal properties.
For instance, for $\alpha=0$, it has been proven
in \cite{CSM} that layer solutions admits nonlocal Modica-type
estimates. Furthermore, in virtue of \cite{cafS},
these equations can be considered as models of a large
variety of nonlocal operators. As a consequence, it is worth studying
the largest possible class of boundary reaction equations. We will
then focus on the following problem:
\begin{equation}\label{eq1-provv} 
\left \{ 
\begin{matrix} 
-{\rm div}\, (\mu(x) \nabla u)+ g(x,u)=0 \qquad 
{\mbox{ on $\R^{n+1}_+$} } \\ 
-\mu(x) u_x=f(u) \qquad 
{\mbox{ on $\partial \R^{n+1}_+$} .}  
\end{matrix} 
\right . 
\end{equation}
under the following structural assumptions (denoted $(S)$):
\begin{itemize}
\item The function $\mu$ is in $L^1((0,r))$, for any $r>0$. Also,
$\mu$ is positive
and bounded over all open sets compactly
contained in $\R^{n+1}_+$, i.e. for all $K \Subset 
\R^{n+1}_+$, 
there exists $\mu_1$,
$\mu_2>0$, possibly depending on $K$,
such that $\mu_1 \leq
\mu(x) \leq \mu_2$, for any $x\in K$.
\item The function~$\mu$ in an $A_2$-Muckenhoupt
weight, that is, there exists $\kappa>0$
such that
\begin{equation}\label{Muck}
\int_a^b \mu(x)\,dx \,
\int_a^b \frac{1}{\mu(x)}\,dx\,\le\,\kappa(b-a)^2 
\end{equation}
for any $b\ge a\ge 0$.
\item The map $(0,+\infty)\ni x\mapsto g(x,0)$
belongs to $L^\infty((0,r))$ for any $r>0$.
Also,
for any $x>0$, the map $\R\ni u\mapsto g(x,u)$
is locally Lipschitz, and given any $R$, $M>0$
there exists $C>0$, possibly depending on $R$ and $M$
in such a way that
\begin{equation}\label{8ikeoqoqoqoo78}
\sup_{{0<x<R}\atop{|u|<M}}|g_u(x,u)|\le C.
\end{equation}
\item The function
$f$ is locally 
Lipschitz in
$\R$.  
\end{itemize} 
In Section \ref{EXT}, using a
Poisson kernel extension,
the fractional equation in~\eqref{lapFrac}
will be reduced to the
extension problem in~\eqref{bdyFrac}, which is
a particular case of \eqref{eq1-provv}.
   
In our setting, \eqref{eq1-provv} may be understood 
in the weak sense, namely supposing that $u\in 
L^\infty_{\rm loc}(\overline{\R^{n+1}_+})$, with  
\begin{equation}\label{hgasj7717177} 
\mu(x)|\nabla u|^2 \in L^1 (B_R^+) 
\end{equation} 
for any $R>0$, 
and that\footnote{Condition \eqref{hgasj7717177} 
is assumed here to make sense of \eqref{eq1}. 
We will see in the forthcoming Lemma \ref{Daf} that it is 
always uniformly fulfilled when $u$ is bounded.

The structural
assumptions on $g$ may be easily
checked when $g(x,u)$ has the product-like
form of $g^{(1)}(x) g^{(2)}(u)$.} 
\begin{equation}\label{eq1} 
\int_{{\R^{n+1}_+}} 
\mu(x)\nabla u\cdot 
\nabla\xi+\int_{{\R^{n+1}_+}} g(x,u)\, \xi= 
\int_{\partial {\R^{n+1}_+}} 
f(u)\xi 
\end{equation} 
for any $\xi:B_R ^+\rightarrow \R$ which is bounded, locally
Lipschitz in the interior of
$\R^{n+1}_+$,
which 
vanishes on $\R^{n+1}_+\setminus B_R$ and such that
\begin{equation}\label{hgasj7717177-bis}
\mu(x)|\nabla\xi|^2
\in L^1 (B_R^+).\end{equation}

As usual,
we are using here the notation~$B_R^+:= B_R 
\cap\R^{n+1}_+$.

In the sequel, we will assume that~$u$ is 
stable, meaning that we suppose 
that 
\begin{equation}\label{sta1} 
\int_{B_R^+} \mu(x)|\nabla\xi|^2+ 
\int_{B_R^+} g_u(x,u)\xi^2 -\int_{\partial B_R^+} 
f'(u)\xi^2\,\ge\,0 
\end{equation} 
for any $\xi$ as above.
 
The stability (sometimes 
also called semistability) condition 
in~\eqref{sta1} appears naturally in the calculus 
of variations setting and it 
is usually related to minimization 
and monotonicity properties. 
In particular,~\eqref{sta1} 
says that the (formal) second variation 
of the energy functional associated 
to the equation has a sign (see, e.g.,~\cite{Moss, FCS, AAC} 
and Section~7 of~\cite{FSV} for further details). 
\bigskip 
 
The main results we prove are a geometric formula, 
of Poincar\'e-type, given in Theorem~\ref{POIN:TH}, 
and a symmetry result, given in Theorem~\ref{SYM:TH}. 
 
For our geometric result, we need to recall 
the following notation. Fixed $x>0$ and~$c\in\R$, we 
look at the level set 
$$ S:= \{  
y\in\R^n {\mbox{ s.t. }} 
u(y,x)=c 
\}.$$ 
We will consider the regular points of~$S$, 
that is, we define 
$$ L:=\{ y\in 
S 
{\mbox{ s.t. }} 
\nabla_y u(y,x)\neq 0 
\}.$$ 
Note that~$L$ depends on the~$x\in(0,+\infty)$ 
that we fixed at the beginning, though we do not keep 
explicit track of this in the notation. 
 
For any point $y\in L$, 
we let $\nabla_L$ to be the tangential gradient 
along~$L$, that is, for any~$y_o\in L$ 
and any~$G:\R^n\rightarrow\R$ smooth in the vicinity of~$y_o$, 
we set 
$$ \nabla_L G(y_o):= 
\nabla_y G(y_o)-\left(\nabla_y G(y_o)\cdot 
\frac{\nabla_y u(y_o,x)}{| 
\nabla_y u(y_o,x)|}\right) 
\frac{\nabla_y u(y_o,x)}{| 
\nabla_y u(y_o,x)|}.$$ 
Since~$L$ is a smooth manifold, in virtue of 
the Implicit Function Theorem (and of the standard 
elliptic 
regularity of $u$ 
apart from the boundary of $\R^{n+1}_+$), 
we can define 
the principal curvatures on it, denoted by 
$$\kappa_1(y,x),\dots, 
\kappa_{n-1}(y,x),$$ for any~$y\in L$. 
We will then define the total curvature 
$$ {\mathcal{K}}(y,x):=\sqrt{ 
\sum_{j=1}^{n-1} \big(\kappa_j (y,x)\big)^2 
}.$$ 
 
We also define 
$$ {\mathcal{R}}^{n+1}_+:=\{ 
(y,x)\in\R^n\times(0,+\infty){\mbox{ s.t. }} 
\nabla_y u(y,x)\neq 0 
\}.$$ 
 
With this notation, we can state 
our geometric formula: 
 
\begin{thm}\label{POIN:TH}
Let $u$ be $C^2_{\rm loc}$ in the interior of $\R^{n+1}_+$.
Assume that 
$u$ is a bounded and stable weak 
solution of \eqref{eq1-provv} under assumptions $(S)$.

Assume furthermore that
for all $r>0$, 
\begin{equation}\label{LipA} |\nabla_y u|\in
{L^\infty(\overline{B_r^+})}. \end{equation}  
Then, 
for any $R>0$ and any~$\phi:\R^{n+1}\rightarrow \R$ which is 
Lipschitz and vanishes on $\R^{n+1}_+\setminus B_R$, we have that 
$$ \int_{ {\mathcal{R}}^{n+1}_+} 
\mu(x) \,\phi^2 \left( 
{\mathcal{K}}^2 |\nabla_y u|^2+ 
\big| 
\nabla_L |\nabla_y u| 
\big| ^2 
\right)\,\leq\, 
\int_{\R^{n+1}_+} 
\mu(x)\, |\nabla_y u|^2 |\nabla \phi|^2 
.$$ 
\end{thm} 

Assumption \eqref{LipA} is natural and it holds in particular in
the important case $g := 0$, $\mu(x)=x^\a$ where $\a \in
(-1,1)$, as discussed in Lemmata~\ref{s8818i1iiii0}
and \ref{ds9882kkk1k1kk1aa}
below.
Interior elliptic regularity also ensures that $u$ is 
smooth inside $\R^{n+1}_+$.

The result in 
Theorem~\ref{POIN:TH} has been inspired 
by the work of~\cite{SZarma, SZcrelle}, as 
developed in~\cite{FAR, FSV}. In particular,~\cite{SZarma, SZcrelle} 
obtained a similar inequality for stable 
solutions of the standard Allen-Cahn equation, and 
symmetry results for possibly singular or degenerate models have 
been obtained in~\cite{FAR, FSV}. 
Related geometric inequalities also played an
important role in \cite{CabCap}.
 
The advantage of the above formula is that 
one  
bounds 
tangential gradients and curvatures of level sets 
of stable solutions in terms of the gradient
of the solution. 
That is, suitable 
geometric quantities of interest 
are controlled by an appropriate 
energy term. 
 
On the other hand, since the geometric formula bounds a 
weighted~$L^2$-norm of any 
test function~$\phi$ by a 
weighted~$L^2$-norm of its 
gradient, we may  
consider Theorem~\ref{POIN:TH} 
as a weighted {P}oincar\'e 
inequality. Again, the advantage of such a formula 
is that the weights have a neat geometric interpretation. 
 
The second result we present is a symmetry result 
in low dimension: 
 
\begin{thm}\label{SYM:TH} 
Let the assumptions of Theorem~\ref{POIN:TH} 
hold and let $n=2$. 

Suppose also
that one of the following conditions~\eqref{g=0}
or~\eqref{g=+} hold, namely assume that
either for any $M>0$ 
\begin{equation}\label{g=0}
{\mbox{the map $(0,+\infty)\ni x\mapsto
\displaystyle\sup_{|u|\le M} |g (x,u)|$ is in $L^1((0,
+\infty))$}} 
\end{equation}
or that
\begin{equation}\label{g=+}
\inf_{{x\in\R^n}\atop{u\in\R}}g(x,u)\,u\,
\ge\,0.
\end{equation}

Suppose also that
there exists $C>0$ in such a way that
\begin{equation}\label{numucr}
\int_0^R \mu(x)\,dx\,\le\, CR^2
\end{equation}
for any $R\ge1$.

Then, there exist~$\omega\in {\rm S}^1$  
and~$u_o: \R\times [0,+\infty)\rightarrow\R$ 
such that 
$$ u(y,x)=u_o(\omega\cdot y,x)$$ 
for any~$(y,x)\in {\R^{3}_+}$. 
\end{thm} 
 
Roughly speaking, Theorem \ref{SYM:TH} 
asserts that, for any $x > 0$, the function $\R^2 
\ni y\mapsto 
u(y,x)$ depends only on one variable. Thus, 
as remarked 
at the beginning of this paper, Theorem \ref{SYM:TH} 
may be seen 
as the analogue of De Giorgi conjecture of~\cite{DeG} 
in dimension $n=2$ 
for equation \eqref{lapFrac}. 

Of course, condition \eqref{numucr} is satisfied,
for instance, for $\mu:=x^\a$ and $\a\in(-1,1)$
and \eqref{g=0}
is fulfilled by $g:=0$, or, more generally,
by $g:=g^{(1)}(x) g^{(2)}(u)$,
with $g^{(1)}$ summable over~$\R^+$
and~$g^{(2)}$ locally Lipschitz.
Also, condition~\eqref{g=+} is fulfilled by~$g:=u^{2\ell+1}$,
with~$\ell\in\N$.

We remark that when $u$ is not bounded,
the claim of Theorem \ref{SYM:TH} 
does not, in general, hold 
(a counterexample being~$\mu:=1$, $f:=0$,
$g:=0$ and~$
u(y_1,y_2,x):=y_1^2-y_2^2$).

Theorem~\ref{aux:P} below
will also provide a result, slightly
more general than Theorem \ref{SYM:TH},
which will be
valid for $n\ge 2$ and without
conditions~\eqref{g=0} or~\eqref{g=+},
under an additional energy assumption.
\bigskip

The pioneering work in~\cite{CSM} 
is related to Theorem~\ref{SYM:TH}. 
Indeed, with different methods,~\cite{CSM} 
proved a result analogous to  
our 
Theorem~\ref{SYM:TH} under the additional 
assumptions that~$\a=0$
and~$f\in C^{1,\beta}$ for some~$\beta>0$ (see, in particular, 
page~1681 and Theorem~1.5 in~\cite{CSM}). 
The method of~\cite{CSM} will be adapted to the case~$\a\in(-1,1)$ 
in the forthcoming paper~\cite{CS}. 

We finally state the symmetry result for equation
\eqref{lapFrac}:
 
\begin{thm}\label{FRAC:TH}
Let $v\in C^2_{\rm loc}(\R^n)$ be a bounded
solution of equation
\eqref{lapFrac}, with $n=2$
and $f$ locally
Lipschitz. 

Suppose that either
\begin{equation}\label{y2v1}
f:=0
\end{equation}
or that
\begin{equation}\label{y2v2}
\partial_{y_2} v>0.
\end{equation}
Then, there exist~$\omega\in {\rm S}^1$  
and~$v_o: \R\rightarrow\R$ 
such that 
$$ v(y)=v_o(\omega\cdot y)$$ 
for any~$y\in \R^2$. 
\end{thm}

The remaining part of the paper is devoted to the proofs 
of Theorems~\ref{POIN:TH} ,~\ref{SYM:TH},~\ref{FRAC:TH}. 
For this, some regularity theory for solutions 
of equation \eqref{eq1-provv} 
will also be needed. 
 
\section{Regularity theory for equation \eqref{eq1-provv}}\label{reg} 
 
This section is devoted to several results we 
need for the regularity theory of equation \eqref{eq1-provv}. 
We do not develop here a complete theory.

We recall that
\begin{equation}
\label{9bis} \mu(x)\,u_x^2
\in L^1(B_R^+)
\end{equation}
for any $R>0$,
due to \eqref{hgasj7717177}.

\subsection{Regularity for equation \eqref{eq1-provv} under 
assumption \eqref{LipA}}

We start with an elementary observation:

\begin{lemma}
If $n=2$ and \eqref{numucr} holds, then
there exists $C>0$ in such a way that
\begin{equation}\label{sii1kkkkj1j1}
\int_{B_{2R}^+ \setminus B_R^+} \mu(x)
\le CR^4
\end{equation}
for any $R\ge1$.
\end{lemma}

\PF Using \eqref{numucr}, we have that
\begin{eqnarray*}
\int_{B_{2R}^+ \setminus B_R^+} \mu(x)
&\le&\int_{0}^{2R} \int_{B_{ 2R} } \mu(x)\,dy\,dx
\\ &\le& C_1 R^2
\int_{0}^{2R} \mu(x)\,dx\\
&\le& C_2 R^4,
\end{eqnarray*}
for suitable $C_1$, $C_2>0$.~\CVD
\medskip

Though not explicitly needed here, we would 
like to point out that the natural integrability 
condition in \eqref{hgasj7717177} 
holds uniformly for bounded solutions.
A byproduct of this gives an
energy estimate, which we will use in the proof of
Theorem \ref{SYM:TH}.
 
\begin{lemma}\label{Daf} 
Let $u$ be a bounded weak 
solution of \eqref{eq1-provv} under
assumptions
$(S)$.

Then, for any $R>0$ there exists $C$, possibly 
depending on $R$, in such a way that 
$$ \| \mu(x) |\nabla u|^2  \|_{L^1 (B_R^+)}\le C.$$ 
Moreover, if 
\begin{itemize}
\item $n=2$,
\item either \eqref{g=0} or~\eqref{g=+} holds,
\item \eqref{numucr}
holds,\end{itemize} then there exists $C_o>0$
such that
\begin{equation}\label{AL}
\int_{B_R^+}\mu(x)\,|\nabla u|^2\,\le\, C_o\, R^2
\end{equation}
for any $R\ge 1$.
\end{lemma} 
 
\PF 
The proof consists just in testing the weak formulation
in~\eqref{eq1}
with 
$\xi:=u
\tau ^2$ where $\tau$ is a cutoff
function
such that $0\le\tau\in C^\infty_0 (B_{2R})$, with $\tau=1$ 
in $B_{R}$ and $|\nabla \tau|\le 8/R$, with $R\ge1$. 
Note that such a $\xi$ is admissible,
since~\eqref{hgasj7717177-bis}
follows from~\eqref{hgasj7717177}.

One then gets from \eqref{eq1} that 
\begin{eqnarray*}&& 
\int_{{\R^{n+1}_+}}
\mu(x)\,
\big( |\nabla u |^2 \tau ^2+2 \tau \nabla u \cdot \nabla \tau 
\big)+ 
\int_{{\R^{n+1}_+}} g(x,u) u \tau^2\\&&\qquad=
\int_{\R^n} f(u) u\tau^2. 
\end{eqnarray*} 
Thus, by Cauchy-Schwarz
inequality,
\begin{eqnarray*}
\int_{\R_+^{n+1}}\mu(x)\,|\nabla u|^2\tau^2
&\le& \frac 12 \int_{\R_+^{n+1}}\mu(x)\,|\nabla u|^2\tau^2
\\&&\quad
+
C_* \Big(
\int_{\R_+^{n+1}}
\mu(x)|\nabla \tau|^2
+\int_{\R^{n}}|f(u)|\,|u|\,\tau^2
\Big)\\&&\quad\quad-
\int_{\R_+^{n+1}}
g(x,u)\,u\,\tau^2,
\end{eqnarray*}
for a suitable constant $C_*>0$.

This, recalling \eqref{g=0}, \eqref{g=+}
and \eqref{sii1kkkkj1j1},
plainly gives the desired result.~\CVD\medskip 
 
We now control further derivatives in $y$, 
thanks to the fact that 
the operator is independent of the variable $y$: 
 
\begin{lemma}\label{Nuovo} 
Let $u$ be a bounded weak 
solution of \eqref{eq1-provv} under conditions~$(S)$.
Suppose that \eqref{LipA} holds.
Then, 
$$\mu(x) |\nabla u_{y_j}|^2 \in L^1(B_R^+)$$ 
for every $R>0$.  
\end{lemma} 
 
\PF 
Given $|\eta|<1$, $\eta\ne 0$, we consider the incremental
quotient
$$ u_\eta(y,x):= \frac{u(y_1,\dots,y_j+\eta,\dots,y_n,x)-
u(y_1,\dots,y_j,\dots,y_n,x)}{\eta}.$$
Since $f$ is locally Lipschitz,
\begin{equation}\label{Al1}
[f(u)]_\eta \le C,
\end{equation}
for some $C>0$, due to \eqref{LipA}.

Analogously, from \eqref{8ikeoqoqoqoo78}
and \eqref{LipA},
for any $R>0$ there exists $C_R>0$ such that
\begin{equation}\label{Al2}
[g(x,u)]_\eta \le C_R
\end{equation}
for any $x\in(0,R)$.

Let now $\xi$ be as requested in~\eqref{eq1}.
Then,~\eqref{eq1} gives that
\begin{eqnarray*}&& 
\int_{{\R^{n+1}_+}}\big[
\mu(x) \nabla u_{\eta}\cdot 
\nabla\xi+\big(g(x,u)\big)_\eta \,\xi\big]- 
\int_{\partial {\R^{n+1}_+}} \big[
f(u)\big]_\eta\xi\\ 
&=& 
-\int_{{\R^{n+1}_+}}\big[
\mu(x) \nabla u\cdot 
\nabla\xi_{-\eta}+g(x,u) \,\xi_{-\eta} 
\big]+
\int_{\partial {\R^{n+1}_+}} 
f(u) \xi_{-\eta}\\ 
&=&0. 
\end{eqnarray*} 
We now consider a smooth cutoff function
$\tau$ such that $0\le\tau\in C^\infty_0 (B_{R+1})$, with $\tau=1$ 
in $B_{R}$ and $|\nabla \tau|\le 2$. 
Taking $\xi:=u_{\eta}\tau^2$ in the above expression, one gets 
\begin{equation}\label{0a8h1hmclakk} 
\begin{split}
& 2 \int_{{\R^{n+1}_+}}
\mu(x) \tau u_{\eta} \nabla u_{\eta
}\cdot \nabla \tau \\&\quad 
+\int_{{\R^{n+1}_+}} 
\mu(x) \tau^2 |\nabla u_{\eta}|^2
+\int_{{\R^{n+1}_+}}  \big(g(x,u)\big)_\eta
u_{\eta}\,\tau^2\\&\quad\quad
=\int_{\partial {\R^{n+1}_+}} 
\big(f(u)\big)_\eta\,u_{\eta} \tau^2.
\end{split}\end{equation}
We remark that the above choice of~$\xi$
is admissible, since
\eqref{hgasj7717177-bis}
follows from~\eqref{LipA} and \eqref{9bis}.
  
Now, by  Cauchy-Schwarz
inequality, we have
\begin{equation*}
\begin{split} 
&\int_{{\R^{n+1}_+}}
\mu(x) \tau u_{\eta} \nabla u_{\eta}\cdot \nabla \tau \geq
-\frac{\varepsilon}{2} \int_{{\R^{n+1}_+}}
\mu(x) \tau^2 |\nabla u_{\eta}|^2\\
& -\frac{1}{2\varepsilon} \int_{{\R^{n+1}_+}}
\mu(x) u_{\eta}^2|\nabla \tau|^2
\end{split}
\end{equation*}
for any $\varepsilon>0$. 

Therefore, by choosing $\varepsilon$ suitably small,~\eqref{0a8h1hmclakk}
reads
\begin{eqnarray*}
&&
\int_{{\R^{n+1}_+}}
\mu(x) \tau^2 |\nabla u_{\eta}|^2
\\
&\le& C\,
\Big[
\int_{B_{R+1}^+}
\mu(x) u_{\eta}^2+
\int_{B_{R+1}^+}  \big| \big(g(x,u)\big)_\eta
u_{\eta}\big|
\\
&&\quad+
\int_{\{|y|\le R\}\times\{x=0\}} \big|
\big(f(u)\big)_\eta
u_{\eta}
\big|\Big].
\end{eqnarray*}
for some $C>0$.

{F}rom~\eqref{LipA}, \eqref{Al1} and \eqref{Al2},
we thus control
$$\int_{B_R^+}
\mu(x) \tau^2 |\nabla u_{\eta}|^2
$$
uniformly in $\eta$.

By sending $\eta\rightarrow 0$
and using Fatou Lemma,
we obtain the desired claim.~\CVD 
\medskip

Following is the regularity needed for
some subsequent computations:

\begin{lemma} 
Let $u$ be $C^2_{\rm loc}$ in the interior of $\R^{n+1}_+$.
Suppose that $u$ is a bounded weak 
solution of \eqref{eq1-provv} under conditions~$(S)$
and that~\eqref{LipA} holds.

Then, 
\begin{equation}\label{SA2} 
\begin{split} 
&{\mbox{for almost any $x>0$, the map $\R^n 
\ni y\mapsto 
\nabla u(y,x)$}}\\ 
&{\mbox{ 
is in $W^{1,1}_{\rm loc}(\R^n, \R^n)$ 
}}\end{split}\end{equation} 
and 
\begin{equation}\label{SA3-provv} 
\begin{split} 
&{\mbox{the map $\R^{n+1}_+ 
\ni (y,x)\mapsto \mu(x) \sum_{j=1}^n 
\big( 
|\nabla u_{y_j}|^2+|u_{y_j}|^2 
\big)$}}\\ 
&{\mbox{is in $L^1(B_r^+)$, for any 
$r>0$. 
}}\end{split}\end{equation} 
What is more,
\begin{equation}\label{SA3}
\begin{split} 
&{\mbox{the map $\R^{n+1}_+ 
\ni (y,x)\mapsto 
\mu(x) \big(  
|\nabla|\nabla_y u||^2+|\nabla_y u|^2 
\big) 
$}}\\ 
&{\mbox{is in
$L^1(B_r^+)$,
for any 
$r>0$. 
}}\end{split}
\end{equation} 
for any $r>0$. 

\end{lemma} 
 
\PF Since $u$ is $C^2_{\rm loc}$
in the interior of $\R^{n+1}_+$,
for any $x\in(\epsilon, 1/\epsilon)$ 
and any $R>0$ 
$$ \int_{B_R} |\nabla u(y,x)| + 
\sum_{j=1}^n |\nabla u_{y_j}(y,x)|\,dy\le C$$ 
for a suitable $C>0$, possibly depending on $\epsilon$ and $R$, 
which proves \eqref{SA2}. 
 
Exploiting Lemma \ref{Nuovo}, \eqref{LipA}
and the local 
integrability of $\mu(x)$, 
one obtains \eqref{SA3-provv}.

To prove \eqref{SA3}, we now
perform the following 
standard approximation argument. 
Define $\Gamma=(\Gamma_1,\dots,\Gamma_n) 
:=\nabla_y u$, 
and let $r$, $\rho>0$ and~$P\in \R^{n+1}_+$ be such 
that~$B_{r+\rho}(P)\subset \R^{n+1}_+$. Fix also $i\in 
\{1,\dots,n+1\}$. 
 
Then, for any $\epsilon>0$, 
\begin{eqnarray*} 
&&\frac{\sum_{j=1}^n \Gamma_j \, 
\partial_i \Gamma_j}{\sqrt{ 
\epsilon^2+ 
\sum_{j=1}^n \Gamma_j^2 
}}\leq 
\frac{2 \, 
|\Gamma| \,|\partial_i \Gamma|}{\epsilon+ 
|\Gamma|}\leq  
2 |\partial_i \Gamma|\in L^1(B_r(P)) 
\\ && 
\lim_{\epsilon\rightarrow 0^+} 
\frac{\sum_{j=1}^n \Gamma_j \,\partial_i \Gamma_j}{\sqrt{ 
\epsilon^2+\sum_{j=1}^n \Gamma_j^2 
}} = \chi_{\{\Gamma\neq 0 \}} 
\frac{ \sum_{j=1}^n \Gamma_j \,\partial_i \Gamma_j }{ 
|\Gamma|} 
\\&&\sqrt{\epsilon^2+\sum_{j=1}^n \Gamma_j^2 
}\leq \epsilon+|\Gamma|\in L^1(B_r(P)) 
\\{\mbox{ and }}\; 
&&\lim_{\epsilon\rightarrow 0^+} 
\sqrt{ 
\epsilon^2+\sum_{j=1}^n \Gamma_j^2 
}=|\Gamma|, 
\end{eqnarray*} 
thanks to \eqref{SA3-provv}. 
As standard, we denote  by $\chi_A$, 
here and in the sequel, the characteristic 
function of a set $A$. 
 
Therefore, by Dominated Convergence Theorem, 
\begin{eqnarray*} 
&& \int_{\R^{n+1}_+}\psi \chi_{\{\Gamma\neq 0 \}} 
\frac{\sum_{j=1}^n \Gamma_j \, 
\partial_i \Gamma_j}{|\Gamma|}= 
\lim_{\epsilon\rightarrow 0^+}\int_{\R^{n+1}_+}\psi 
\frac{\sum_{j=1}^n \Gamma_j \, 
\partial_i \Gamma_j}{\sqrt{ 
\epsilon^2+\sum_{j=1}^n \Gamma_j^2}}\\ 
&&\quad= 
\lim_{\epsilon\rightarrow 0^+} 
\int_{\R^{n+1}_+} 
\psi \,\partial_i \left( \sqrt{ 
\epsilon^2+\sum_{j=1}^n \Gamma_j^2 
}\right)\\&&\quad\quad 
=- 
\lim_{\epsilon\rightarrow 0^+} 
\int_{\R^{n+1}_+} 
(\partial_i\psi)\sqrt{ 
\epsilon^2+\sum_{j=1}^n \Gamma_j^2 
}\\&&\quad\quad\quad=- 
\int_{\R^{n+1}_+} 
(\partial_i\psi) |\Gamma|. 
\end{eqnarray*} 
for any $\psi\in C^\infty_0 (B_r(P))$. 
 
Thus, 
since~$P$, $r$ and~$\rho$ can be arbitrarily chosen, 
we have that 
$$ \partial_i |\Gamma|= 
\chi_{ \{\Gamma\neq 0 \} } 
\frac{\sum_{j=1}^n \Gamma_j  
\partial_i\, \Gamma_j}{|\Gamma| 
}$$ 
weakly and almost everywhere in $\R^{n+1}_+$. 
 
Accordingly, 
\begin{eqnarray*} 
&& |\nabla |\nabla_y u||^2=|\nabla |\Gamma||^2 
=\sum_{i=1}^{n+1}(\partial_i |\Gamma|)^2\\ 
&&\quad 
\leq 
\sum_{i=1}^{n+1} 
\left( 
\frac{\sum_{j=1}^n \Gamma_j  
\,\partial_i \Gamma_j}{|\Gamma| 
} 
\right)^2\leq 
\sum_{i=1}^{n+1} 
|\partial_i \Gamma |^2\\ 
&&\quad\quad= 
\sum_{i=1}^{n+1} 
\sum_{j=1}^{n} 
(\partial_i u_{y_j})^2= 
\sum_{j=1}^{n} 
|\nabla u_{y_j}|^2.\end{eqnarray*} 
Then, \eqref{SA3-provv} implies
\eqref{SA3}.~\CVD

\subsection{Verification of assumption \eqref{LipA}}
In this section, we 
show that \eqref{LipA} is always
satisfied in
the important case $g :=0$, $\mu(x):=x^\a$,
with $\a\in(-1,1)$.

More precisely, we state the following result,
the proof of which can be found in \cite{CS}:

\begin{lemma}\label{lhah} 
Let $u$ be a bounded weak solution
of \eqref{bdyFrac} and 
assume that $f$ is locally Lipschitz. 
Then there exists a constant $C>0$ 
depending on $R$ and $\beta \in (0,1)$ such that  
\begin{itemize} 
\item the function $u$ is 
H\"older-continuous of exponent $\beta$ and 
$$\|u\|_{C^\beta(\overline{B_{R}^+})} \leq C,$$ 
\item for all $j=1,...,n$, the function $u_{y_j}$ 
is H\"older-continuous of exponent $\beta$ and 
\begin{equation}\label{s8828181aa} 
\|u_{y_j}\|_{C^\beta(\overline{B_{R}^+})} \leq C. 
\end{equation} 
\end{itemize} 
\end{lemma} 
 
We can now prove the following gradient bound,
which says that \eqref{LipA} holds for bounded solutions
of equation \eqref{bdyFrac}:
  
\begin{lemma}\label{s8818i1iiii0} 
Let $u$ be a bounded weak 
solution of \eqref{bdyFrac} and assume 
that $f$ is locally Lipschitz.  
 
Then, there exists a constant $C>0$ such that  
$$\|\nabla_y u \|_{L^\infty
(\overline{\R^{n+1}_+})} \leq C.$$ 
\end{lemma} 
 
\PF {F}rom \eqref{s8828181aa}, 
$\nabla_y u $ is bounded 
in, say, $\overline{\R^{n+1}_+} \cap \left \{ 0 \leq x \leq 3\right \}$. 
 
Now, in $\overline{\R^{n+1}_+} \cap \left \{ x > 3\right \}$, 
equation \eqref{bdyFrac} is nondegenerate and 
therefore, the gradient bound follows from standard elliptic 
theory.~\CVD 
   
\section{Proof of Theorem \ref{POIN:TH}} 
 
Besides few technicalities, the 
proof of Theorem \ref{POIN:TH} consists simply 
in plugging the right test function in stability 
condition~\eqref{sta1} and in using the linearization 
of \eqref{eq1-provv} to get rid of the unpleasant terms. 
Following are the rigorous details of the proof. 
   
By~\eqref{SA2}, we have that 
$$ \int_{\R^{n+1}_+} \mu(x) \nabla u_{y_j} 
\cdot \Psi = 
\int_0^\infty \mu(x) \int_{\R^n} \nabla u_{y_j} 
\cdot \Psi\,dy\,dx= 
-\int_{\R^{n+1}_+}  \mu(x) \nabla u 
\cdot \Psi_{y_j}$$ 
for any~$j=1,\dots, n$ and any~$\Psi\in C^\infty 
(\R^{n+1}_+, \R^n)$ supported in~$B_R$. 
 
Thus,
making use of~\eqref{eq1},
we conclude that 
\begin{equation}\label{a711aa} 
\begin{split} 
&\int_{\R^{n+1}_+} \mu(x) \nabla u_{y_j} 
\cdot \nabla\psi\\ 
&\quad=- 
\int_{\R^{n+1}_+} \mu(x) \nabla u 
\cdot \nabla\psi_{y_j}\\ &\quad\quad 
=-\int_{\partial \R^{n+1}_+} f(u)\psi_{y_j}+\int_{\R^{n+1}_+} 
g(x,u)  \psi_{y_j} \\ 
&\quad\quad\quad 
=\int_{\partial \R^{n+1}_+}
\big(f(u)\big)_{y_j}\psi-\int_{\R^{n+1}_+} g_u(x,u)u_{y_j} \psi \\ 
&\quad\quad\quad\quad 
=\int_{\partial \R^{n+1}_+} f'(u)u_{y_j} \psi-\int_{\R^{n+1}_+}
g_u(x,u) u_{y_j} \psi 
\end{split} 
\end{equation} 
for any~$j=1,\dots, n$ and any~$\psi\in 
C^\infty 
(\R^{n+1}_+)$ supported in~$B_R$. 
 
A density argument (recall~\eqref{Muck}
and see, e.g.,
Lemma~3.4, Theorem~2.4 and~(2.9) in~\cite{CPSC}),
via 
%% Lemma \ref{s8818i1iiii0} and
\eqref{SA3-provv}, implies that~\eqref{a711aa} 
holds for~$\psi:=u_{y_j} \phi^2$, 
where~$\phi$ is 
as in the statement of 
Theorem~\ref{POIN:TH}, therefore 
\begin{equation}\label{ssj9818919} 
\begin{split} 
&\int_{\partial B_R^+} f'(u)|\nabla_y u|^2 \phi^2 \\
&\quad= 
\sum_{j=1}^n\int_{B_R^+} \mu(x) \nabla u_{y_j}\cdot 
\nabla (u_{y_j}\phi^2)+\sum_{j=1}^n\int_{B^+_R} g_u(x,u) u^2_{y_j}\phi^2 \\ 
&\quad= 
\sum_{j=1}^n\int_{B_R^+} \mu(x) \left(|\nabla u_{y_j}|^2\phi^2+ 
u_{y_j}\nabla u_{y_j}\cdot\nabla\phi^2\right)\\&
\qquad+\sum_{j=1}^n\int_{B^+_R}
g_u(x,u) u^2_{y_j}\phi^2 \\ 
&\quad=\int_{B_R^+} 
\mu(x) \left(\sum_{j=1}^n 
|\nabla u_{y_j}|^2\phi^2+ 
\phi\nabla\phi\cdot\nabla|\nabla_y
u|^2\right)\\&
\qquad+
\int_{B^+_R} g_u(x,u) |\nabla_y u|^2
\phi^2 . 
\end{split} 
\end{equation} 
Now, we make use of~\eqref{sta1} 
by taking~$\xi:=|\nabla_y u|\phi$ 
(this choice was also performed 
in~\cite{SZarma, SZcrelle, FAR, FSV}; 
note that \eqref{LipA} 
and~\eqref{SA3} imply \eqref{hgasj7717177-bis}
and so they
make it possible to 
use here such a test function). We thus obtain 
\begin{eqnarray*} 0&\leq& \int_{B_R^+} 
\mu(x) \Big( \big| \nabla |\nabla_y u|\big|^2 \phi^2+ 
|\nabla_y u|^2 |\nabla\phi|^2\\
&&+2|\nabla_y u|\phi\nabla\phi\cdot 
\nabla|\nabla_y u| \Big)\\ 
&& +\int_{B_R^+}g_u(x,u) |\nabla_y u|^2 \phi^2-\int_{\partial B_R^+} 
f'(u)|\nabla_y u|^2\phi^2. 
\end{eqnarray*} 
This and~\eqref{ssj9818919} 
imply that 
\begin{equation}\begin{split} 
\label{s88818181} 
&\int_{\R^{n+1}_+} \mu(x)\phi^2 \Big( 
\sum_{j=1}^{n-1} \big(\partial_x u_{y_j}\big)^2 
- 
\big(\partial_x |\nabla_y u| 
\big)^2 \\&
+\sum_{j=1}^{n-1} |\nabla_y u_{y_j}|^2-\big| 
\nabla_y |\nabla_y u| 
\big|^2 
\Big)\\ 
&+\int_{\R^{n+1}_+} \mu(x) 
\phi\nabla\phi\cdot\left( 
\nabla|\nabla_y u|^2 
- 
2|\nabla_y u| 
\nabla|\nabla_y u|\right) 
\\&\qquad= 
\int_{B_R^+} \mu(x)\phi^2 \Big( 
\sum_{j=1}^{n-1} |\nabla u_{y_j}|^2-\big| 
\nabla|\nabla_y u| 
\big|^2 
\Big)\\ 
&\qquad 
+\int_{\R^{n+1}_+} \mu(x)\phi\nabla\phi\cdot 
\left(\nabla|\nabla_y u|^2 
- 
2|\nabla_y u| 
\nabla|\nabla_y u|\right) 
\\&\qquad\qquad 
\le \int_{\R^{n+1}_+} 
\mu(x) |\nabla\phi|^2|\nabla_y u|^2 
.\end{split}\end{equation} 
 
Let now~$r$, $\rho>0$ and~$P\in \R^{n+1}_+$ be such 
that~$B_{r+\rho}(P)\subset \R^{n+1}_+$. 
We consider~$\gamma$ to be either~$|\nabla_y u|$ 
or~$u_{y_j}$. In force of~\eqref{SA3-provv}
and~\eqref{SA3}, we see that~$\gamma$ 
is in~$W^{1,2}( B_r(P))$, and so in~$W^{1,1}_{\rm loc} 
(B_r (P))$. 
 
Thus, by 
Stampacchia Theorem (see, e.g., Theorem~6.19 
in~\cite{LOSS}), 
$\nabla \gamma=0$ for almost 
any~$(y,x)\in B_r(P)$ such that~$\gamma(y)=0$. 
 
Hence, since~$P$, $r$ and~$\rho$ can be 
chosen arbitrarily, 
we have that 
$\nabla |\nabla_y u| =0=\nabla u_{y_j}$ 
for almost every~$(y,x)$ such that~$\nabla_y u(y,x)=0$. 
 
Accordingly,~\eqref{s88818181} 
may be written as 
\begin{eqnarray*} 
&& \int_{{{\mathcal{R}}^{n+1}_+}} \mu(x)\phi^2 \left( 
\sum_{j=1}^{n-1} \big(\partial_x u_{y_j}\big)^2 
- 
\big(\partial_x |\nabla_y u| 
\big)^2\right)\\ 
&&+ 
\int_{{\mathcal{R}}^{n+1}_+} 
\mu(x)\phi^2\left( 
\sum_{j=1}^{n-1} |\nabla_y u_{y_j}|^2-\big| 
\nabla_y|\nabla_y u| 
\big|^2 
\right) 
\\ 
&\le& \int_{\R^{n+1}_+} 
\mu(x) |\nabla\phi|^2|\nabla_y u|^2 
. 
\end{eqnarray*} 
Therefore, by standard differential geometry formulas 
(see, for example, 
equation~(2.10) in~\cite{FSV}), we obtain 
\begin{eqnarray} 
\label{yuiooaoo} 
\nonumber&& 
\int_{{{\mathcal{R}}^{n+1}_+}} \mu(x)\phi^2 \left( 
\sum_{j=1}^{n-1} \big(\partial_x u_{y_j}\big)^2 
- 
\big(\partial_x |\nabla_y u| 
\big)^2\right)\\ 
&&+ 
\int_{ {\mathcal{R}}^{n+1}_+} 
\mu(x) \phi^2 \left( 
{\mathcal{K}}^2 |\nabla_y u|^2+ 
\big| 
\nabla_L |\nabla_y u| 
\big| ^2 
\right) 
\\ 
&\le& \int_{\R^{n+1}_+} 
\mu(x) |\nabla\phi|^2|\nabla_y u|^2 
.\nonumber 
\end{eqnarray} 
We now note that, on~${{\mathcal{R}}^{n+1}_+}$, 
$$ \big(\partial_x |\nabla_y u| 
\big)^2 =\left| 
\frac{\nabla_y u \cdot \nabla_y u_x}{\nabla_y u} 
\right|^2\leq |\nabla_y u_x|^2= 
\sum_{j=1}^{n-1} \big(\partial_x u_{y_j}\big)^2. 
$$ 
This and~\eqref{yuiooaoo} 
complete the proof of 
Theorem~\ref{POIN:TH}.~\CVD 
 
\section{Proof of Theorem \ref{SYM:TH}} 
 
The strategy for proving Theorem \ref{SYM:TH} 
is to test the geometric formula of 
Theorem~\ref{POIN:TH} against an appropriate capacity-type 
function to make the left  
hand side vanish. 
This would give that the curvature of the level sets for fixed~$x>0$ 
vanishes and so that these level sets are flat, as desired 
(for this, 
the vanishing of the tangential gradient term is also 
useful to take care of the possible 
plateaus of~$u$, where 
the level sets are not smooth manifold: see Section~2.4 in~\cite{FSV} 
for further considerations). 
 
Some preparation is needed for the proof 
of Theorem \ref{SYM:TH}. 
Indeed, Theorem \ref{SYM:TH} will 
follow from the subsequent Theorem~\ref{aux:P}, which 
is valid for any dimension~$n$ and without the restriction
in either \eqref{g=0}
or \eqref{g=+}. 
 
We will use the notation~$X:=(y,x)$ for points in~$\R^{n+1}$.  
 
Given~$\rho_1\le\rho_2$, we also define 
$${\mathcal{A}}_{\rho_1,\rho_2}:=\{ 
X\in\R^{n+1}_+{\mbox{ s.t. }}|X|\in [\rho_1,\rho_2] 
\}.$$ 
 
\begin{lemma}\label{tatay} 
Let~$R>0$ 
and~$h:B_R^+\rightarrow\R$ be a nonnegative 
measurable function.  
 
For any~$\rho\in (0,R)$,
let 
$$ \eta(\rho):=\int_{B^+_{\rho}} h.$$ 
Then, 
$$\int_{{\mathcal{A}}_{\sqrt R, R}}\frac{h(X)}{|X|^2}\,dX 
\leq 2\int_{\sqrt R}^R t^{-3}\eta(t)\,dt+\frac{\eta(R)}{R^2}. 
$$ 
\end{lemma} 
 
\PF  
By Fubini Theorem, 
\begin{eqnarray*} 
&&\int_{{\mathcal{A}}_{\sqrt R, R}}\frac{h(X)}{2 |X|^2}\,dX 
\\&=& 
\int_{{\mathcal{A}}_{\sqrt R, R}}\int_{|X|}^R 
t^{-3} h(X)\,dt\,dX 
+\int_{{\mathcal{A}}_{\sqrt R, R}}\frac{h(X)}{2 R^2}\,dX 
\\&=&\int_{\sqrt R}^R 
\int_{{\mathcal{A}}_{\sqrt R, t}} 
t^{-3} h(X)\,dX\,dt+\frac{1}{2R^2} 
\int_{{\mathcal{A}}_{\sqrt R, R}} {h(X)}\,dX 
\\ &\leq& 
\int_{\sqrt R}^R 
\int_{B^+_t} 
t^{-3} h(X)\,dX\,dt+\frac{1}{2R^2} 
\int_{B^+_R}{h(X)}\,dX 
, 
\end{eqnarray*} 
from which we obtain the desired result.~\CVD 
 
\begin{thm}\label{aux:P} 
Let $u$ be as requested in Theorem \ref{POIN:TH}.
Assume furthermore that  
there exists~$C_o\geq 1$ in such a way that 
\begin{equation}\label{en:bound} 
\int_{B^+_R} \mu(x) |\nabla u|^2\le C_o\,  
R^2\end{equation} 
for any~$R\ge C_o$. 
 
Then there exist~$\omega\in {\rm S}^{n-1}$ 
and~$u_o: \R\times(0,+\infty)\rightarrow\R$ 
such that 
$$ u(y,x)=u_o(\omega\cdot y,x)$$ 
for any~$(y,x)\in\R^{n+1}_+$. 
\end{thm} 
 
\PF {F}rom Lemma~\ref{tatay} (applied here with
$h(X):=\mu(x) |\nabla 
u(X)|^2$) and \eqref{en:bound}, we obtain
\begin{equation}\label{7s77s88} 
\int_{{\mathcal{A}}_{\sqrt R, R}}\frac{\mu(x) |\nabla u(X)|^2 
}{ |X|^2}\leq C_1\log R 
\end{equation} 
for a suitable~$C_1$, 
as long as~$R$ is large enough. 
 
Now we define 
$$ \phi_R(X):=\left\{ 
\begin{matrix} 
\log R & {\mbox{ if $|X|\le \sqrt R$,}}\\ 
2\log\big( R/|X|\big)\Big)  
& {\mbox{ if $\sqrt R<|X|< R$,}} 
\\ 
0 & {\mbox{ if $|X|\ge R$}} 
\end{matrix} 
\right.$$ 
and we observe that 
$$|\nabla\phi_R|\leq \frac{C_2\,\chi_{ 
{\mathcal{A}}_{\sqrt R, R} 
}}{|X|},$$ 
for a suitable~$C_2>0$. 
 
Thus, plugging~$\phi_R$ inside the geometric 
inequality of Theorem~\ref{POIN:TH}, we obtain 
\begin{eqnarray*} 
&& (\log R)^2\int_{B^+_{\sqrt{R}}\cap 
{\mathcal{R}}^{n+1}_+ 
} 
\mu(x) \left( 
{\mathcal{K}}^2 |\nabla_y u|^2+ 
\big| 
\nabla_L |\nabla_y u| 
\big| ^2 
\right)\\&&\qquad\qquad\,\leq\,C_3 
\int_{ 
{\mathcal{A}}_{\sqrt R, R} 
}\frac{ 
\mu(x) |\nabla_y u|^2}{|X|^2} 
\end{eqnarray*} 
for large~$R$. 
 
Dividing by~$(\log R)^2$, 
employing~\eqref{7s77s88} 
and taking~$R$ arbitrarily large, we see that 
$$ 
{\mathcal{K}}^2 |\nabla_y u|^2+ 
\big| 
\nabla_L |\nabla_y u| 
\big| ^2 
$$ 
vanishes identically on~${\mathcal{R}}^{n+1}_+$, 
that is~${\mathcal{K}}=0=\big| 
\nabla_L |\nabla_y u| 
\big|$ on~${\mathcal{R}}^{n+1}_+$. 
 
Then, the desired result follows 
by Lemma~2.11 of~\cite{FSV} (applied to the function~$y\mapsto 
u(y,x)$, for any fixed~$x>0$).~\CVD 
\medskip 

We now complete the proof of 
Theorem~\ref{SYM:TH}.
We observe that, under the assumptions of
Theorem~\ref{SYM:TH},
estimate \eqref{en:bound}
holds, thanks to \eqref{AL}.
Consequently, the hypotheses of Theorem~\ref{SYM:TH} 
imply the ones of Theorem~\ref{aux:P}, 
from which the claim in Theorem~\ref{SYM:TH} follows.~\CVD 

\section{Proof of Theorem \ref{FRAC:TH}}\label{EXT}

We use Theorem \ref{SYM:TH} to prove Theorem
\ref{FRAC:TH}. For this, given a function
$v$ satisfying \eqref{lapFrac}, 
we select an extension\footnote{The extension
is not, in general, unique. For instance,
both the functions $u:=0$ and
$u:=x^{1-\a}$
satisfy $\,{\rm div}\,(x^\a \nabla u)=0$ in $\R^{n+1}_+$
with $u=0$ on $\partial\R^{n+1}_+$.}
$u$ satisfying \eqref{bdyFrac} by a suitable Poisson
kernel, whose theory has been developed in \cite{cafS}.

For this, we use the following result of~\cite{cafS}:
\begin{lemma}
The function 
$$P(y,x)=C_{n,\a}\frac{x^{1-\a}}{\Big ( x^2+|y|^2
\Big)^\frac{n+1-\a}{2}}$$
is a solution of 
\begin{equation}\label{POI}
\left\{ 
\begin{matrix} 
-{\rm div}\, (x^\a \nabla P)=0 \qquad 
{\mbox{ on $\R^n\times(0,+\infty)$}} 
\\ 
P=\delta_0 
\qquad{\mbox{ on $\R^n\times\{0\}$,}}\end{matrix} 
\right.\end{equation} 
where $\a \in (-1,1)$ and $C_{n,\a}$ is a normalizing constant such
that 
$$\int_{R^n} P(y,x)\,dy=1. $$
\end{lemma} 

We now come to the proof of Theorem \ref{FRAC:TH}. Let $v$ be a
bounded solution of \eqref{lapFrac} and consider the function 
\begin{equation}\label{conv}
u(y,x)= \int_{\R^n}P(y-z,x)\,v(z)\,dz=
\int_{\R^n}P(\xi,x)\,v(y-\xi)\,d\xi. 
\end{equation} 
Note that since $P(x,.)\in L^1(\R^n)$
and $v \in L^\infty(\R^n)$ and
by the embedding $L^1\,*\,L^\infty
\subset L^\infty$, we have that $u$
is bounded in $\R^{n+1}_+$
if $v$ is bounded in $\R^n$. 

We now prove the following regularity result. 

\begin{lemma}\label{ds9882kkk1k1kk1aa}
Let $v$ be
bounded and
$C^2_{\rm loc}(\R^n)$.
Let $u$ be given by \eqref{conv}.

Then, for all $R>0$ there exists a constant
$C>0$ such that 
$$\|x^\a u_x \|_{L^\infty(\overline{B_R^+})} \leq C. $$
\end{lemma}

\PF
Since $P$ has unit mass, we have the relation 
$$ u(y,x) - v(y) =
C_{n,\a} \int_{\R^n} \frac{x^{1-\a}
(v(y-\xi)-v(y))
}{( x^2 +|\xi|^2)^{(n+1-\a)/2}} d\xi. $$
Therefore,
\begin{equation*}
\begin{split}
&x^\a u_x = x^\a \partial_x (u(y,x)-v(y))
\\
&\qquad=C_{n,\a} \int_{\R^n} \frac{ [ (1-\a) |\xi|^2 -n x^2]\,
(v(y-\xi)-v(y))
}{( x^2 +|\xi|^2)^{(n+3-\a)/2}} d\xi 
\end{split}
\end{equation*}
This bounds the quantity $x^\a u_x$ by
$$ \int_{\R^n} \frac{
|v(y-\xi)-v(y)|
}{( x^2 +|\xi|^2)^{(n+1-\a)/2}} d\xi $$
which
is controlled by
\begin{eqnarray*} &&\int_{\R^n} \frac{
|v(y-\xi)-v(y)|
}{|\xi|^{(n+1-\a)}} d\xi\\&&\qquad \leq
\int_{ |\xi|\ge 1} \frac{
2\|v\|_{L^\infty(\R^n)}
}{|\xi|^{(n+1-\a)}} d\xi
+
\int_{ |\xi|\le 1} \frac{
\|\nabla v\|_{L^\infty( B_1 (y) )}
}{|\xi|^{(n-\a)}} d\xi.
\end{eqnarray*}

The last two terms are summable and one gets  the bound
$$ \| x^\a u_x\|_{L^\infty (\overline{B_R^+})}
\le C \big(\|v\|_{L^\infty(\R^n)} +
\|\nabla v\|_{L^\infty( B_{R+1})}\big),$$  
as desired.~\CVD
\medskip

We now complete the proof of Theorem \ref{FRAC:TH}
via the following argument.
We take $u$ as defined in \eqref{conv}
and we observe that \eqref{LipA} and \eqref{9bis}
are satisfied, thanks to the local integrability
of $x^{-\a}$ and Lemmata~\ref{s8818i1iiii0}
and~\ref{ds9882kkk1k1kk1aa}.

Also, $u$ is stable, because of
either~\eqref{y2v1}
or~\eqref{y2v2}. 

Indeed, if~\eqref{y2v1}
holds, then~\eqref{sta1}
is obvious since~$f:=0=:g$
in this case. 

If, on the other hand,~\eqref{y2v2}
holds, then~$u_{y_2}= P*v_{y_2}>0$ in~$\R^{n+1}_+$, 
and~$u_{y_2}(y,0)=v_{y_2}(y)
>0$ on~$\partial\R^{n+1}_+$, thanks to Lemma~\ref{lhah}.

Therefore,
given~$\xi:B_R ^+\rightarrow \R$ which is bounded, locally
Lipschitz in the interior of
$\R^{n+1}_+$,
which
vanishes on
$\R^{n+1}_+\setminus B_R$ and such
that~\eqref{hgasj7717177-bis}
holds, we
use~\eqref{a711aa}
with~$\psi:=\xi^2/u_{y_2}$ (here, 
$j:=2$, $g:=0$, $\mu:=x^\a$
and~\eqref{SA3-provv} make the choice of
such a $\psi$
admissible), and we get
$$ \int_{\partial\R^{n+1}_+} f'(u)\xi^2=
\int_{\R^{n+1}_+}
2x^\a\xi\frac{\nabla u_{y_2}\cdot\nabla\xi}{u_{y_2}}
- x^\a \xi^2 \frac{|\nabla u_{y_2}|^2}{u_{y_2}^2}
.$$
This, by Cauchy-Schwarz
inequality, gives~\eqref{sta1}
and so $u$
is stable.\medskip

Then,
we apply Theorem \ref{SYM:TH} to get that~$u(y,x)=
u_o(\omega\cdot y,x)$ for any $y\in\R^2$
and any $x>0$,
for an appropriate direction~$\omega$.

By Lemma~\ref{lhah}, $u$ is continuous up to~$\{x=0\}$
and so~$u(y,0)=
u_o(\omega\cdot y,0)$.

Since, by \eqref{POI} and \eqref{conv},
$$u|_{\partial \R^{n+1}_+}=v,$$
the proof of
Theorem~\ref{FRAC:TH}
is complete.~\CVD

\section*{Acknowledgments} 
 
This collaboration has started on the occasion of 
a very pleasant visit of the authors at the 
University of Texas at Austin. We would like to thank the 
department of mathematics for its kind hospitality.  
 
EV has been partially supported by~{\em MIUR 
Me\-to\-di va\-ria\-zio\-na\-li ed equa\-zio\-ni 
dif\-fe\-ren\-zia\-li non\-li\-nea\-ri}.

\bigskip\bigskip\bigskip 
 
{\small Preprint~\cite{FSV} is available on line at  
 
{\tt 
http://www.math.utexas.edu/mp\_arc/}  
 
no. 07-171.} 
 
\bigskip\bigskip\bigskip\bigskip 
 
{\em YS} --  
Universit\'e Aix-Marseille 3, Paul C\'ezanne -- 
LATP -- 
Marseille, France. 
 
{\tt sire@cmi.univ-mrs.fr} 
\medskip 
 
{\em EV} -- 
Universit\`a di Roma Tor Vergata -- 
Dipartimento di Matematica -- 
I-00133 Rome, Italy. 
 
{\tt valdinoci@mat.uniroma2.it}

\begin{thebibliography}{CPSC94}

\bibitem[AAC01]{AAC}
Giovanni Alberti, Luigi Ambrosio, and Xavier Cabr{\'e}.
\newblock On a long-standing conjecture of {E}. {D}e {G}iorgi: symmetry in 3{D}
  for general nonlinearities and a local minimality property.
\newblock {\em Acta Appl. Math.}, 65(1-3):9--33, 2001.
\newblock Special issue dedicated to Antonio Avantaggiati on the occasion of
  his 70th birthday.

\bibitem[ABS98]{ABS}
Giovanni Alberti, Guy Bouchitt{\'e}, and Pierre Seppecher.
\newblock Phase transition with the line-tension effect.
\newblock {\em Arch. Rational Mech. Anal.}, 144(1):1--46, 1998.

\bibitem[AC00]{AC}
Luigi Ambrosio and Xavier Cabr{\'e}.
\newblock Entire solutions of semilinear elliptic equations in {$ \mathbb R\sp
  3$} and a conjecture of {D}e {G}iorgi.
\newblock {\em J. Amer. Math. Soc.}, 13(4):725--739 (electronic), 2000.

\bibitem[BCN97]{BCN}
H.~Berestycki, L.~A. Caffarelli, and L.~Nirenberg.
\newblock Monotonicity for elliptic equations in unbounded {L}ipschitz domains.
\newblock {\em Comm. Pure Appl. Math.}, 50(11):1089--1111, 1997.

\bibitem[Ber96]{B}
Jean Bertoin.
\newblock {\em L\'evy processes}, volume 121 of {\em Cambridge Tracts in
  Mathematics}.
\newblock Cambridge University Press, Cambridge, 1996.

\bibitem[CC06]{CabCap}
Xavier Cabr{\'e} and Antonio Capella.
\newblock Regularity of radial minimizers and extremal solutions of semilinear
  elliptic equations.
\newblock {\em J. Funct. Anal.}, 238(2):709--733, 2006.

\bibitem[CCFS98]{POW}
M.~Chipot, M.~Chleb{\'{\i}}k, M.~Fila, and I.~Shafrir.
\newblock Existence of positive solutions of a semilinear elliptic equation in
  {${\mathbb R}\sp n\sb {+}$} with a nonlinear boundary condition.
\newblock {\em J. Math. Anal. Appl.}, 223(2):429--471, 1998.

\bibitem[CG08]{Mar}
Luis Caffarelli and Mar{\'{\i}}a del~Mar Gonz{\'a}lez.
\newblock {G}amma convergence of an energy functional related to the fractional
  {L}aplacian.
\newblock {\em Preprint}, 2008.

\bibitem[CPSC94]{CPSC}
Valeria Chiad{\`o}~Piat and Francesco Serra~Cassano.
\newblock Relaxation of degenerate variational integrals.
\newblock {\em Nonlinear Anal.}, 22(4):409--424, 1994.

\bibitem[CRS07]{CRS}
Luis Caffarelli, Jean-Michel Roquejoffre, and Yannick Sire.
\newblock Free boundaries with fractinal {L}aplacians.
\newblock {\em In preparation}, 2007.

\bibitem[CS07a]{CS}
X.~Cabr\'e and Y.~Sire.
\newblock Semilinear equations with fractional {L}aplacians.
\newblock {\em In preparation}, 2007.

\bibitem[CS07b]{cafS}
Luis Caffarelli and Luis Silvestre.
\newblock An extension problem related to the fractional laplacian.
\newblock {\em Commun. in PDE}, 32(8):1245, 2007.

\bibitem[CSM05]{CSM}
Xavier Cabr{\'e} and Joan Sol{\`a}-Morales.
\newblock Layer solutions in a half-space for boundary reactions.
\newblock {\em Comm. Pure Appl. Math.}, 58(12):1678--1732, 2005.

\bibitem[CSS08]{CSS}
Luis Caffarelli, Sandro Salsa, and Luis Silvestre.
\newblock {R}egularity estimates for the solution and the free boundary of the
  obstacle problem for the fractional {L}aplacian.
\newblock {\em To appear in Invent. Math.}, 2008.

\bibitem[CT04]{CT}
Rama Cont and Peter Tankov.
\newblock {\em Financial modelling with jump processes}.
\newblock Chapman \& Hall/CRC Financial Mathematics Series. Chapman \&
  Hall/CRC, Boca Raton, FL, 2004.

\bibitem[CV06]{CV}
Luis Caffarelli and Alexis Vasseur.
\newblock Drift diffusion equations with fractional diffusion and the
  quasi-geostrophic equation.
\newblock {\em Preprint}, 2006.

\bibitem[DG79]{DeG}
Ennio De~Giorgi.
\newblock Convergence problems for functionals and operators.
\newblock In {\em Proceedings of the International Meeting on Recent Methods in
  Nonlinear Analysis (Rome, 1978)}, pages 131--188, Bologna, 1979. Pitagora.

\bibitem[DL76]{DL}
G.~Duvaut and J.-L. Lions.
\newblock {\em Inequalities in mechanics and physics}.
\newblock Springer-Verlag, Berlin, 1976.
\newblock Translated from the French by C. W. John, Grundlehren der
  Mathematischen Wissenschaften, 219.

\bibitem[Far02]{FAR}
Alberto Farina.
\newblock Propri\'et\'es qualitatives de solutions d'\'equations et syst\`emes
  d'\'equations non-lin\'eaires.
\newblock 2002.
\newblock Habilitation \`a diriger des recherches, Paris VI.

\bibitem[FCS80]{FCS}
Doris Fischer-Colbrie and Richard Schoen.
\newblock The structure of complete stable minimal surfaces in {$3$}-manifolds
  of nonnegative scalar curvature.
\newblock {\em Comm. Pure Appl. Math.}, 33(2):199--211, 1980.

\bibitem[FJK82]{FJK}
E.~Fabes, D.~Jerison, and C.~Kenig.
\newblock The {W}iener test for degenerate elliptic equations.
\newblock {\em Ann. Inst. Fourier (Grenoble)}, 32(3):vi, 151--182, 1982.

\bibitem[FKS82]{FKS}
Eugene~B. Fabes, Carlos~E. Kenig, and Raul~P. Serapioni.
\newblock The local regularity of solutions of degenerate elliptic equations.
\newblock {\em Comm. Partial Differential Equations}, 7(1):77--116, 1982.

\bibitem[FSV07]{FSV}
Alberto Farina, Bernadiro Sciunzi, and Enrico Valdinoci.
\newblock Bernstein and {D}e giorgi type problems: new results via a geometric
  approach.
\newblock {\em Preprint}, 2007.

\bibitem[GG98]{GG1}
N.~Ghoussoub and C.~Gui.
\newblock On a conjecture of {D}e {G}iorgi and some related problems.
\newblock {\em Math. Ann.}, 311(3):481--491, 1998.

\bibitem[GG03]{GG2}
Nassif Ghoussoub and Changfeng Gui.
\newblock On {D}e {G}iorgi's conjecture in dimensions 4 and 5.
\newblock {\em Ann. of Math. (2)}, 157(1):313--334, 2003.

\bibitem[Giu84]{giusti}
Enrico Giusti.
\newblock {\em Minimal surfaces and functions of bounded variation}, volume~80
  of {\em Monographs in Mathematics}.
\newblock Birkh\"auser Verlag, Basel, 1984.

\bibitem[Lan72]{landkof}
N.~S. Landkof.
\newblock {\em Foundations of modern potential theory}.
\newblock Springer-Verlag, New York, 1972.
\newblock Translated from the Russian by A. P. Doohovskoy, Die Grundlehren der
  mathematischen Wissenschaften, Band 180.

\bibitem[LL97]{LOSS}
Elliott~H. Lieb and Michael Loss.
\newblock {\em Analysis}, volume~14 of {\em Graduate Studies in Mathematics}.
\newblock American Mathematical Society, Providence, RI, 1997.

\bibitem[MP78]{Moss}
William~F. Moss and John Piepenbrink.
\newblock Positive solutions of elliptic equations.
\newblock {\em Pacific J. Math.}, 75(1):219--226, 1978.

\bibitem[Muc72]{muck}
Benjamin Muckenhoupt.
\newblock Weighted norm inequalities for the {H}ardy maximal function.
\newblock {\em Trans. Amer. Math. Soc.}, 165:207--226, 1972.

\bibitem[Sav08]{savin}
Ovidiu Savin.
\newblock Phase transitions: Regularity of flat level sets.
\newblock {\em To appear in Ann. of Math.}, 2008.

\bibitem[SZ98a]{SZarma}
Peter Sternberg and Kevin Zumbrun.
\newblock Connectivity of phase boundaries in strictly convex domains.
\newblock {\em Arch. Rational Mech. Anal.}, 141(4):375--400, 1998.

\bibitem[SZ98b]{SZcrelle}
Peter Sternberg and Kevin Zumbrun.
\newblock A {P}oincar\'e inequality with applications to volume-constrained
  area-minimizing surfaces.
\newblock {\em J. Reine Angew. Math.}, 503:63--85, 1998.

\end{thebibliography}
\end{document}